\documentclass[12pt]{article}
\usepackage[english]{babel}
\usepackage[utf8]{inputenc}
\usepackage[T1]{fontenc}

\usepackage{amsfonts, amsmath, natbib}
\bibpunct{[}{]}{;}{a}{}{,~}
\usepackage[hidelinks]{hyperref}

\usepackage[left=1.0 in,right=1.0 in,marginparwidth=0cm,
   top=1.5 cm,height=238mm,headheight=15pt]{geometry}

\DeclareMathOperator{\Ee}{\mathrm{E}}

\DeclareMathOperator{\Pp}{\mathrm{P}}

\title{Early Highlights in the History of the Bernstein-von Mises Theorem}
\date{2025\\ November}
\author{Hans Fischer\thanks{Formerly: Katholische Universität Eichstätt-Ingolstadt, Bavaria. Mail: hans.fischer@ku.de}}
\begin{document}
\maketitle
\thispagestyle{empty}

\noindent The designation ``Bernstein-von Mises theorem'' (abbreviated as BvM theorem) is apparently due to Lucien Le Cam, as will be explained below more closely. Roughly, the assertion of this theorem is that the posterior distribution of a parameter, conditioned on a large sample, is approximately normal, independent of a particular prior.

\section{The Basic Idea}

\label{basicidea}

The basic idea behind the BvM theorem can be described for one-dimensional parameters as follows: Let $X_1, X_2, \dots, X_n$ be a sample of i.i.d. r.v's with values in $\mathbb{R}$, each of which obeys a probability density $f_\theta$ with respect to a measure $\nu$ on $\cal{B}(\mathbb{R})$, where $f_\theta(x) = f(\theta, x)$ depends on the one-dimensional parameter $\theta \in (\theta_1, \theta_2)$. Let $\pi: (\theta_1, \theta_2) \to \mathbb{R}_0^+$ be the a-priori density (with respect to the Lebesgue-measure) of $\theta$.

For $\theta_1 \leq z_1 < z_2 \leq \theta_2$ (provided the respective integrals exist) we have:
\begin{equation} \label{basic}
\Pp(\theta \in (z_1, z_2)| X_1=x_1, \dots, X_n=x_n) = \frac{\int_{z_1}^{z_2} f_\theta(x_1) \cdots f_\theta(x_n) \pi(\theta) d \theta}{\int_{\theta_1}^{\theta_2} f_\theta(x_1) \cdots f_\theta(x_n) \pi(\theta) d \theta}
\end{equation} if the denominator on the right side is $\ne 0$.
Let $$p_n(\theta | x_1, \dots, x_n):= f_\theta(x_1) \cdots f_\theta(x_n) \pi(\theta)$$ and
$$\ell(\theta | x_1, \dots, x_n) := \sum \log(f_\theta(x_i)).$$ We assume that there exists a parameter $\theta_n \in (\theta_1, \theta_2)$ such that $\ell'(\theta_n | x_1, \dots, x_n) =0$ and $\ell''(\theta_n | x_1, \dots, x_n) = - \alpha_n^2 < 0$ (derivatives with respect to $\theta$). $\theta_n$ is a maximum likelihood estimate for $\theta$, and $- \alpha_n^2$ as a sum of $n$ terms is of order $n$ under certain conditions. Therefore, if $f(\theta, x)$ is sufficiently smooth with respect to $\theta$, we obtain an expansion, for $x/\alpha_n$ in a neighborhood of $\theta_n$:
\begin{equation} \label{expans} \ell(\theta_n + \frac{x}{\alpha_n}) = \ell(\theta_n) - \frac{x^2}{2} + R_n(x).\end{equation} If we can show that $R_n(x)$ becomes sufficiently ``small'' for large $n$, then, under the additional supposition of a continuous $\pi$, we obtain :
\begin{equation} \label{rough} p_n(\theta_n + \frac{x}{\alpha_n} | x_1, \dots, x_n) \approx f_{\theta_n}(x_1) \cdots f_{\theta_n}(x_n) \pi(\theta_n) e^{- \frac{x^2}{2}}.\end{equation}  This is already a primitive statement on the asymptotic normality of the posterior.

In order to arrive at a proper limit theorem for the posterior, particular assumptions about the behavior of the maximum likelihood estimate $\theta_n$ and about $f_\theta(x)$ are necessary. In the simplest case we suppose that $\theta_n \to \theta_0 \in (\theta_1, \theta_2)$ for $n \to \infty$ while $\pi(\theta)$ is continuous and positive for $\theta= \theta_0$. If we are able to show that $\exp(R_n(x))$ tends to 1 for all relevant $x$ and this limit is compatible with integration, and if we can infer $\alpha_n = \mbox{O}(\sqrt{n})$, then, according to (\ref{basic}):
\begin{multline} \label{limit} \Pp(\theta \in (\theta_n+\frac{a}{\alpha_n}, \theta_n+\frac{b}{\alpha_n})| X_1=x_1, \dots, X_n=x_n) =
\frac{\int_a^b p_n(\theta_n + \frac{x}{\alpha_n}| x_1, \dots, x_n) dx}{\int_{(\theta_1-\theta_n)\alpha_n}^{(\theta_2-\theta_n)\alpha_n} p_n(\theta_n + \frac{x}{\alpha_n} | x_1, \dots, x_n) dx} \\
\to \frac{\pi(\theta_0) \int_a^b \exp(-\frac{x^2}{2}) dx}{\pi(\theta_0) \int_{-\infty}^\infty \exp(-\frac{x^2}{2}) dx} = \frac{1}{\sqrt{2 \pi}}  \int_a^b \exp(-\frac{x^2}{2}) dx. \end{multline}

Already in this rather simple setting there are many ``ifs'' behind the limit assertion. In some special cases of sample distributions, as binomial or multinomial distributions, only assumptions about the limit behavior of $\theta_n$ are essentially needed.  In the discussion of the general case, which was started by Le Cam in 1953, more intricate considerations are required, also because of the peculiarities of non-elementary conditional probabilities and of the different possible modes of convergence.

\section{Laplace's Approximations}

\label{lapapprox}

The present paper is devoted to contributions up to Le Cam's 1953 paper, and all these contributions were in the tradition of Laplace's pertinent work. This work has been discussed in detail by several historians (see, e.g., \citep{hald:1998, dale:1999, hald:2007, gorro:2016}). Therefore it should be sufficient to recapitulate at this place only a few of its most significant ideas.

Already in \citeyearpar{laplace:1774b} Laplace (1749--1827) published one of his most prominent papers on probability, in which he also derived an approximate normal posterior in the binomial case. The parameter is the success probability $\theta \in (0,1)$, the $X_i$ are two-valued with 1 for ``success'' and 0 for ``failure''. As always, Laplace applies a uniform prior only. Then  relation (\ref{basic}) becomes:
\begin{equation} \label{binomcase} \Pp(\theta \in (z_1, z_2)| X_1=x_1, \dots, X_n=x_n) = \frac{\int_{z_1}^{z_2} x^s (1-x)^{n-s} dx}{\int_{0}^{1} x^s (1-x)^{n-s} dx}, \mbox{ where } s=\sum x_i.\end{equation} From a modern point of view, the relative frequency $\theta_n= s/n$ would be the maximum likelihood estimate for $\theta$. Following the Bernoullian law of large numbers, for Laplace it is a matter of course, however, to take the relative frequency as an estimator for $\theta$, and he shows that $x=s/n$ maximizes the integrands. Moreover, $\alpha_n$ is equal to $\sqrt{\frac{n}{\theta_n(1-\theta_n)}}$ in the binomial case. Laplace applies an expansion very similar to (\ref{expans}), but he does not explicitly rescale the argument by means of $\alpha_n$ (this was done in contributions from 1900 on only). Thus, he arrives at a result equivalent to
$$\ell(\theta_n + z) = \ell(\theta_n) - \frac{\alpha_n^2 z^2}{2} + R_n(z),$$ where $R_n(z)$ is indicated by the first term of a series with the order of magnitude $\mbox{O}(n z^3)$. Under the assumption $z^3 = \mbox{o}(1/n)$ the approximation (according to (\ref{rough}))
$$p_n(\theta_n + z | s) \approx \left(\frac{s}{n}\right)^s \left(\frac{n-s}{n}\right)^{n-s} e^{- \frac{n^3 z^2}{2s(n-s)}}$$ is followed. Contrary to the procedure described above, the integral $\int_{0}^{1} x^s (1-x)^{n-s} dx$ is evaluated by elementary methods, and the result $\frac{(n+1)!}{s!(n-s)!}$ is approximated by Stirling's formula. Altogether, Laplace obtains for $\omega^3 = \mbox{o}(1/n)$ an approximation equivalent to:
$$\Pp(\theta_n - \omega \leq \theta \leq \theta_n + \omega | s) \approx \frac{1}{\sqrt{2\pi}} \int_{-\omega \alpha_n}^{\omega \alpha_n} \exp(-\frac{x^2}{2}) dx.$$

In his 1774 publication, Laplace was primarily interested in the case that $\omega \alpha_n$ despite the smallness of $\omega$ could indefinitely grow if $n$ grows. Based on a result of Euler he showed that
$\int_{-\infty}^\infty \exp(-x^2) dx = \sqrt{\pi},$ and therefore, for large $n$ and arbitrarily small $\varepsilon$:
$$\Pp(\theta_n - \varepsilon \leq \theta \leq \theta_n + \varepsilon | s) \approx 1.$$ Whereas Jakob Bernoulli's law of large numbers proved the closeness of relative frequency and success probability for the situation of an experiment which has still to be executed, Laplace had thus established a corresponding result for the a posteriori situation in which the outcomes of the experiment are known. Later on, this result would be named ``second law of large numbers'' (this designation was applied by \citet{mises:1919a, mises:1919b}, for example).

Laplace's further work on inverse probability in which normal posteriors occurred, dealt beside the binomial case (that was applied to various statistical problems) with the parameters of multinomial distributions and location parameters of continuous distributions. The most significant sources are \citep{laplace:1781, laplace:1785} and Laplace's main work, the {\it Théorie analytique des probabilités} \citeyearpar{laplace:1812, laplace:1820a}. In these sources also general analytic explanations are given about approximating integrals of the type
$$\int_c^d (f_1(x))^{s_1} \cdots (f_n(x))^{s_n} g(x) dx$$ with non-negative functions and ``large'' exponents $s_i$ such that the whole integrand has a unique maximum in $(c,d)$.  Also corresponding integrals in two variables are treated. These integrals are very similar to the integrals in (\ref{basic}). $g(x)$ may even be interpreted as an a-priori density from today's point of view. The difference to the procedure indicated in our section (1) is, however, that there the maximum refers to the product of the $f$'s only. And, actually, Laplace in the context of probabilistic problems proper never discussed a non-uniform prior.

In his discussion of the multinomial model Laplace restricted himself to three-valued r.v's. He derived an approximation of the posterior by a bivariate normal distribution, but he did not give an application of his formula to a statistical problem.

The normal posterior for location parameters was only treated in the {\it Théorie analytique} (2nd book, No. 23), in the context of a problem which had repeatedly been discussed by Laplace since his earliest papers: to find,  in some sense, an optimal estimate for a quantity $a$ on the basis of direct observations $x_i$ which are flawed by errors $\epsilon_i$, such that $x_i = a + \epsilon_i$. For the i.i.d. observational errors Laplace assumed a density $f(x)$ proportional to $\exp(-\psi(x^2))$, where $\psi$ is monotonically increasing (including the case of a constant $\psi$). Then (again with a uniform prior)
$p_n(a | x_1, \dots, x_n) = \prod f(x_i-a)$, and the density $\tilde p_n$ of the posterior is proportional to $p_n$. Laplace applied the condition
$$\int_{-\infty}^\infty |\hat a - a| \tilde p_n(a | \overline x_1, \dots, \overline x_k) = \mbox{min}$$ for the estimate of $a$. This was not new, he had this condition already in his celebrated 1774 paper. The condition leads to the median of $\tilde p_n$ (also shown in 1774 already), and, because of the particular shape of $f$, to its mode. Laplace derived the asymptotic normality for $\tilde p_n$, but the median of the normal distribution depended in an unknown way on $\psi$, and therefore he could not find a general formula or method for computing $\hat a$ as a function of the $x_i$ by these a posteriori considerations. Only in the special case of $f$ being a normal density with zero expectation it could be shown that this median is equal to the arithmetic mean of the $x_i$.\\

\noindent For a better understanding of Laplace's work the following remarks may be helpful:\\

1.) The now usual standards of ``classical'' analysis became only common from ca. 1870 on. Instead of discussing the limit behavior of residual terms, as usual now, Laplace restricted himself to the observation that the first terms of a series expansion expressing the approximation error become ``very small''  when the number of sample elements $n$ becomes ``very large.''

2.) Accordingly, Laplace did not derive limit theorems, as in section (\ref{basicidea}),  but approximations for ``large numbers.'' For such limit theorems centering and re-scaling of the pertinent random variables is typical (in our case by means of $\theta_n$ and $\alpha_n$). A first limit theorem of this type occurs in Chebyshev's work on the central limit theorem (1887). \citet{markov:1900} used the corresponding norming in the binomial case when deriving approximations to the posterior under the assumption of a uniform prior.

3.)  Laplace's restriction to uniform priors can perhaps be explained as follows. In the {\it Essai philosophique}, his work on philosophical aspects of probability, which also served as an introduction to his {\it Théorie analytique} from the third edn. on (1820), Laplace also comprehensively expounded his theory of the association of ideas. And in this context he explained how to obtain probabilities empirically. If you want to explore the probability of an event which can be repeated indefinitely often, and if you start in a state of total ignorance, this corresponds to a uniform prior. And when you repeat the trials, your knowledge becomes better and better, a process which is modeled by the ``second law'' of large numbers. And if you start in a state of an already partial knowledge, this state has been preceded by a process starting from total ignorance. Therefore it should be sufficient to consider uniform priors.

The further development during the 19th century and the first decade of the 20th century did not yield significant progress---from a purely mathematical point of view. A certain focus was on extending Laplace's results to multi-dimensional parameters (Bienaymé 1838, Karl Pearson and Filon 1898 for multinomial models; Edgeworth 1908/09 for location parameters, see \citep[chapt. 9--10]{hald:2007}). Analytical rigor was not essential in these papers, however.

For statistics, especially Bienaymé's results are remarkable, because he succeeded, in a certain sense, in deriving a result which Laplace had failed to reach: that the arithmetic mean is the asymptotic median of the posterior of the location parmeter. Admittedly, Bienaymé considered a somewhat different situation: The observational errors were assumed to be restricted to discrete values $\gamma_1, \dots, \gamma_n$, and under this assumption it was possible to base the solution of the problem on the normal approximation of a multinomial distribution.

 The paradigm of uniform prior was only occasionally weakened. At least
\citet[\S 95]{cournot:1843} gave a plausibility argument for the fact that in the binomial case under the assumption of a large number of trials the influence of a non-uniform prior was restricted to a small interval around the observed relative frequency, with the consequence that the prior---approximately---could be considered as constant.
(Hald's \citeyearpar[chapt. 11.5]{hald:2007} reference to \citet{kries:1886} in this context is not correct, a seemingly similar remark in \citep[125]{kries:1886} concerns a fundamentally different situation in which one single observation is subjected to an error law $\varphi(x)$ that quickly decreases outside of $x=0$.)

\section{Bernstein's Rigorous Proof for the Binomial}

Sergei Natanovich Bernshtein (1880--1968) was born in Odessa. He studied in Göttingen and Paris, where he obtained a PhD in 1904 with a thesis on Hilbert's 19th problem. From 1907 on he taught at the university of Kharkov (Kharkiv), at the beginning with some difficulties, because his academic degree earned in Paris was not acknowledged in Russia. From about 1910 on, Bernstein also devoted himself to probability,  apparently being strongly influenced by Chebyshev's, Lyapunov's, and Markov's work. Due to his affinity to French language, Bernshtein in many cases used the French transcript of his name: ``Bernstein.''

One of Bernstein's most prominent students was Jerzy Neyman (1894--1981), who  as an undergraduate studied in Kharkiv from 1914 to 1917 (due to problems with his eyes Neyman was not admitted to the army in the first world war).
Regarding  Bernstein's lectures on probability that Neyman attended, we find an interesting remark in Reid's Neyman biography \citeyearpar[24]{reid:1982} which resulted from many meetings and conversations of the author with Neyman:
\begin{quote} \normalsize
Many years later, Neyman, by then identifying himself quite completely
with the frequentist point of view, described ``a particular attack on the Bayes
front for which I personally have great respect. .. [which] I learned of in
Bernstein's lectures in 1915 or 1916.'' Here he was referring to a theorem of
Bernstein's (to come up later in his own work) which provides interesting
information on the behavior of the Bayesian and frequentist approaches in
situations where large numbers of observations are made.)
In addition to attending Bernstein's lectures on probability theory,
Neyman--as a member of the Presidium of the Mathematics Club--also
helped to prepare them for mimeographed publication.
\end{quote}
The ``particular attack'' will be discussed below. The ``theorem of
Bernstein's'' was a special version of the BvM theorem.

Lecture notes on Bernstein's probability course appeared in 1917 in a lithographed version with title ``Kurs leksii po teorii veroyatnostei,'' but without reference to Neyman. A student named Zshukovok is mentioned in the title page, and the subtitle is (in English translation) ``edition of the association for mutual support of the students of mathematics at Kharkov university.'' Apparently, these lecture notes remained widely unknown outside Russia. They must not be considered as a precursor of Bernstein's textbook ``Teoriya veroyatnostei'' (first edn. 1927, fourth edn. 1946), which expounded the theory in a rather elementary way only.

In some points the lecture notes reflect a very up-to-date state of the art in the years after 1910: Bernstein's own axiomatics, the central limit theorem in Lyapunov's version, and a rigorous and very general treatment of the BvM-theorem in the binomial case. At least this version was referred to in ``Teoriya veroyatnostei'' (first edn. 1927, footnote on p. 271), and a complete proof (very similar to the original proof in the lecture notes, as we will see) can finally be found in the fourth edn. of this textbook (1946), in appendix 4 (reprinted in {\it Sobranie sochinenii} IV, pp. 448--452). Apparently, all these sources remained widely unknown outside Russia, until Le Cam in the bibliography of his 1953 article hinted at the 1917 work, if in a somewhat incomplete and therefore confusing manner. We will come back to this issue in the final section of the present paper.

In the 1917 lithograph, Bernstein's derivation of the BvMt theorem is on pages 118--122. Detailed assumptions are introduced ad hoc in his proof. The assertion (with our abbreviations introduced in sect. \ref{basicidea} and sect. \ref{lapapprox}) is:
$$\Pp\left(x_1 \beta_n < \theta - \theta_n < x_2 \beta_n | s \right) \to \frac{1}{\sqrt{\pi}} \int_{x_1}^{x_2} \exp(-x^2) dx, \quad \beta_n := \frac{\sqrt{2}}{\alpha_n}.$$
Bernstein starts from the relation
$$\Pp\left(z_1 < \theta < z_2 | s \right) = \frac{\int_{z_1}^{z_2} x^s (1-x)^{n-s} f(x) dx}{\int_{0}^{1} x^s (1-x)^{n-s} f(x) dx},$$
 corresponding to (\ref{binomcase}), but now with a general prior $f$, and first arrives at
$$\Pp(\theta \in (\theta_n+ x_1 \beta_n, \theta_n+x_2 \beta_n)| s) =
\frac{\int_{y_1}^{y_2} p_n(\theta_n + y| s) dy}{\int_{-\theta_n}^{1-\theta_n} p_n(\theta_n + y | s) dy}, $$
where $y_{1/2}:= x_{1/2}\beta_n$, and
$$p_n(\theta_n + y| s) := \theta_n^s (1-\theta_n)^{n-s} f(\theta_n+y) (1+ y/\theta_n)^s(1-y/(1-\theta_n))^{n-s}.$$ This corresponds to the first line of (\ref{limit}).

At his place Bernstein introduces his main assumptions: $\theta_n = s/n$  tends to a limit $\theta$ in the open interval $(0,1)$ as $n$ tends to infinity, and $f$ is assumed to have the property \begin{equation} \label{f-prop} f(\theta_n + z_n) = f(\theta_n)(1+\alpha')\end{equation} such that, if $z_n \to 0$ for $n \to \infty$, the quantity $\alpha'$ tends to 0 also. This holds if $f$ is continuous in $\theta$ and $>0$ there.

Now, the denominator is split off into three integrals with the domains $(-\theta_n, -\varepsilon)$, $(-\varepsilon, \varepsilon)$, $(\varepsilon, 1-\theta_n)$, and the mean value theorem for integrals is applied, with the consequence
\begin{equation} \label{bruch} \Pp(\theta \in (\theta_n+ y_1, \theta_n+y_2)| s) = \frac{f(\theta_n + \overline y) \int_{y_1}^{y_2} (1+ y/\theta_n)^s(1-y/(1-\theta_n))^{n-s} dy}{f(\theta_n+y') \int_{-\varepsilon}^{\varepsilon} (1+ y/\theta_n)^s(1-y/(1-\theta_n))^{n-s} dy + R},\end{equation}  where $R$ is the sum of the two remaining integrals;  $\overline y$ is between $y_1$ and $y_2$, and $y'$ is between $-\varepsilon$ and $\varepsilon$.

By aid of Taylor's theorem, Bernstein obtains for $|y| \leq \varepsilon$:
\begin{equation} \label{expansion} \log\left((1+ y/\theta_n)^s(1-y/(1-\theta_n))^{n-s}\right) = - \frac{n^3y^2}{2 s(n-s)}(1+\alpha_1),\end{equation} where $\alpha_1$ is a function of $y$ and $\alpha_1 = \mbox{O}(\varepsilon)$. Decisive for this latter property is that $\theta_n$ remains bounded away from 0 and 1 for sufficiently large $n$ because its limit $\theta$ is not equal to 0 or 1.

Bernstein assumes $\varepsilon = L\beta_n$ with $L = n^{1/12}$. Then we have $\varepsilon = \mbox{O}(n^{-5/12})$. Other choices would have been possible such that $L \to \infty$ and $L\beta_n \to 0$. Bernstein even admits that the limits $x_1, x_2$ in the assertion of the theorem are between $-L$ and $L$ (and not necessarily fixed).

An estimate for $R$ is obtained as follows: We have
$$R = \int_{-s/n}^{-\varepsilon} f(\frac{s}{n}+y) r_n(y) dy + \int_{\varepsilon}^{(1-s/n)} f(\frac{s}{n}+y) r_n(y) dy, \quad r_n(y) := (1+\frac{ny}{s})^s (1- \frac{ny}{n-s})^{n-s}.$$ $r_n(y)$ is motonotically increasing from $y=-s/n$ to $y=0$, and decreasing from $y=0$ to $y=(n-s)/n$. Therefore, on account of (\ref{expansion}), outside the interval $(-\varepsilon, \varepsilon)$ the inequality $$r_n(y) \leq \exp\left(- \frac{n^3\varepsilon^2}{2 s(n-s)}(1+\alpha_1)\right)$$ is valid. Bernstein now further assumes $f(z) \le M$ for all $z$, and thus arrives at
\begin{equation} \label{R-absch} |R| < M \int_{-s/n}^{1-s/n} \exp\left(- \frac{n^3\varepsilon^2}{2 s(n-s)}(1+\alpha_1)\right) dy = M \exp(-L^2(1+\alpha_1)).\end{equation}

On account of (\ref{f-prop}) and (\ref{expansion}) the right-hand side of (\ref{bruch}) becomes
\begin{equation} \label{bruch1} \frac{(1+\alpha') \int_{x_1}^{x_2} \exp(-x^2 (1+\alpha_1)) dx}{(1+\alpha'') \int_{-L}^{L} \exp(-x^2 (1+\alpha_1)) dx + R/(\beta_n f(\theta_n))}.\end{equation} Because of (\ref{R-absch}) $R/(\beta_n f(\theta_n))$ tends to 0. Again by aid of the mean value theorem for integrals and because of $\alpha', \alpha'' \to 0$, Bernstein finally shows that (\ref{bruch1}) tends to
$$ \frac{\int_{x_1}^{x_2} \exp(-x^2) dx}{\int_{-\infty}^{\infty} \exp(-x^2) dx},$$ which completes the proof.

To some extent, a more streamlined version of the proof can be found in the fourth edn. of Bernstein's textbook, as already hinted at. Here the author does not assume $\theta_n \to \theta$, but he presupposes that $a < \theta_n <1-a$ for a positive $a$, and that the prior $f(x)$ is continuous in $[0, 1]$ and has a positive lower bound in $(a, 1-a)$.

It would have been possible to significantly simplify the proof by use of Lebesgue's theorem on dominated convergence. Of course, Bernstein was aware of this theorem in 1917. It is to be suspected that he did not apply it for didactical reasons.

\section{Von Mises's 1919 Paper}

Originally educated as an engineer, Richard von Mises (1883--1953) had become more and more interested in the mathematical side of this discipline, and from 1919 on he held a professorship at the University of Straßburg (the capital of Alsace, until 1919 under German administration)  for applied mathematics. Despite his duties in the World War---he was significantly involved in the organization of the air force of the Austrian Empire---toward the end of the war he submitted two large papers to {\it Mathematische Zeitschrift}, which were devoted to probability, a field he had only occasionally dealt with until this time. The first of these papers \citep{mises:1919a} was on ``fundamental theorems of probability theory,'' the second \citep{mises:1919b} on ``basic notions of probability theory.'' Both articles would become very influential on the further development of mathematical probability. While by the second von Mises became the leading proponent of frequentism, the first determined methods and style of many later works by various authors, most notably by his way of applying ``distribution functions'' and Stieltjes integrals, and his use of characteristic functions (named ``komplexe Adjunkte'' by him).

In order to reach a joint theory on which all theorems concerning normal limit distributions could be based, von Mises did something that was very smart: In a first part of the ``fundamental theorems'' he established auxiliary theorems on the convergence of products of functions to a function of the type $e^{-x^2}$.

The two most important of these auxiliary theorems are as follows: In the one-dimensional case the product $p_n(u)$ is considered, where
$$p_n(u):= f_1(a_1 + \frac{u}{r_n}) \cdots f_n(a_n + \frac{u}{r_n}).$$
$a_i$ is a sequence of real numbers, the (sufficiently smooth) functions $f_i$ are defined on $\mathbb{R}$, $f_i(a_i) = 1$, $f_i'(a_i) = 0$, $f_i''(a_i) =: -2s_i^2 < -s^2 < 0$, and $\quad r^2_n:= \sum_{i=1}^n s_i^2$. The functions $f_i$ are subject to some further conditions regarding boundedness and integrability. If an additional function $\psi(x)$ is bounded and continuous in $x=a$, then for all $-\infty \leq x_1 < x_2 \leq +\infty$:
\begin{equation} \label{vMtheorem1}
\int_{x_1}^{x_2} \psi(a + \frac{u}{r_n}) f_1(a_1 + \frac{u}{r_n}) \cdots f_n(a_n + \frac{u}{r_n}) du \to \psi(a) \int_{x_1}^{x_2} e^{-u^2} du \quad (n \to \infty).\end{equation}

The multi-dimensional case could in principle be treated like the one-dimensional. The $f_i$ are now defined in $\mathbb{R}^t$, and for the vectors $a_k \in \mathbb{R}^t$ it is assumed that $f(a_k) = 1$, $\nabla f(a_k) = 0$,
$\frac{\partial^2}{\partial x_i \partial x_j} f_k(a_k) =: -2s_k^{i,j}$ where $\sum_{i,j =1}^t s_k^{i,j}y_i y_j \geq 0$ for all vectors $y \in \mathbb{R}^t$. Further properties of the $f_k$ and the single function $\psi$ are presupposed in analogy to the one-dimensional case.

In order to simplify the exposition, von Mises restricted the assertion of the multi-dimensional theorem to the case where the functions obey the condition that (in modern notation) the matrix $H_n$ with matrix elements $h_n^{i,j}:= \frac{1}{n} \sum_{k=1}^n s_k^{i,j}$ remains constant for all $n$. Von Mises gave the sketch of a proof---he indicated the main arguments in analogy to the one-dimensional case---that
\begin{equation} \label{vMtheorem2}
\int_{Z} \psi(a+\frac{z}{\sqrt{n}}) f_1(a_1+\frac{z}{\sqrt{n}}) \cdots f_n(a_n + \frac{z}{\sqrt{n}}) dz \to \psi(a) \int_Z
\exp(- z^T H_n z) dz, \end{equation} where $Z$ is a ``finite or infinite part of space.''

Von Mises had very good knowledge and skills in classic, ``post-Weierstrassian'' analysis, but he was apparently not informed about more recent theories of measure and integration (Lebesgue, Radon). Therefore many of his proofs were rather tedious, and it would have been possible to simplify them and to reach more general results. As his references show, he had some (superficial) knowledge of Laplace's {\it Théorie analytique}, and he knew some works of Chebyshev and Markov, especially the latter's textbook on probability, whose second edition had been translated into German \citeyearpar{markov:1912}. Von Mises further referred to Bienaymés important 1838 paper and to Bruns's book of 1906 \citep{bienayme:1838, bruns:1906}. Yet, he was apparently not acquainted with the newest results of probability theory as one can see from the fact that he did not refer to Lyapunov's papers \citeyearpar{lyapunov:1900, lyapunov:1901} on the central limit theorem.

A considerable part of the ``fundamental theorems'' was dedicated to central limit theorems, though. In this part von Mises reached significant innovations with respect to local limit theorems for random variables both with continuous and with lattice distributions, by applying his auxiliary theorems to characteristic functions. In another part, limit theorems on inverse probabilities were discussed, and in this field we are mainly interested here.

The binomial case of the BvM theorem was mastered by von Mises in a way which was, in principle, equivalent to Bernstein's, but differed regarding several details, especially the use of auxiliary theorems. The starting point was the problem of drawing black and white balls from an urn with replacement. The probability for a white ball was designated by $x$ with $0<x<1$, and it was assumed that in $n$ drawings $an$ white balls have occurred, where $0<a<1$ (of course, von Mises assumed $[an]$ white balls, but we will use his notation in the following). $a$ was presupposed to be independent of $n$, and this was a certain restriction. Following Laplace, in a first step the prior was assumed to be uniform. Then the probability density $w_n(x)$ of the posterior  was
$$w_n(x) = C_n[x^a(1-x)^{1-a}]^n,$$ where $C_n$ was determined by the
condition
\begin{equation} \label{normierung} \int_0^1 w_n(x) dx = 1.\end{equation} Von Mises applied the one-dimensional auxiliary theorem to functions
$$f_k(x) = \left\{\begin{array}{crl} \left(\frac{x}{a}\right)^a \left(\frac{1-x}{1-a}\right)^{1-a} & \mbox{ for } & 0 < x < 1\\
0 & \mbox{else.} & \end{array} \right.
$$ Then $f_k$ follow the conditions of the auxiliary theorem for $a_k = a$ and $r_n^2=n/(2a(1-a))$, and therefore we have
\begin{equation} \label{miseslimit} \frac{1}{C_n a^{an} (1-a)^{(1-a)n}} \int_{x_1}^{x_2} w_n(a + \frac{u}{r_n}) du \to \int_{x_1}^{x_2} e^{-u^2} du.\end{equation} If we set $w_n(a + \frac{u}{r_n})=0$ for $a+u/r_n \not \in [0,1]$, then because of (\ref{normierung}) the equality $$r_n \int_{-\infty}^\infty w_n(a + \frac{u}{r_n}) du = 1$$ ensues. Therefore,  from (\ref{miseslimit})  with $x_1 = -\infty$, $x_2 = \infty$ the limit
$$\frac{1}{r_nC_n a^{an} (1-a)^{(1-a)n}} \to \sqrt{\pi}$$ can be followed. The consequence is the limit statement for the a-posteriori probability of the success probability $x$:
\begin{equation} \label{vMbinom} \Pp(a-\frac{x_1}{r_n} < x < a+\frac{x_1}{r_n} | an) = r_n \int_{x_1}^{x_2}  w_n(a + \frac{u}{r_n}) du \to \frac{1}{\sqrt{\pi}} \int_{x_1}^{x_2} e^{-u^2} du.\end{equation}

In a second step von Mises admitted a general a-priori density $\psi(x)$ which had to be continuous for all $x$ and positive in $x=a$ (it would have been sufficient to assume it to be bounded and continuous in $x=a$, and $\psi(a) > 0$). Now he was able to simply recall the general version of his auxiliary theorem in writing:
\begin{quote} \dots it follows that $\lim w_n$ deviates from the previous value only by the constant factor $\psi(a)$. Because this factor, which we suppose to be different from zero, in  the determination of $C_n$ cancels out, we see that the results [essentially (\ref{vMbinom})] remain unchanged \dots
\end{quote}

The multinomial case with $t$ different events was tackled by von Mises in an analogous way. Here the posterior density for the probabilities $x_1, \dots, x_{t-1}$ and $x_t = 1-\sum_{i=1}^{t-1} x_i$ under the supposition of constant relative frequencies $a_1, \dots, a_{t-1}$ and $a_t = 1-\sum_{i=1}^{t-1} a_i$ (all probabilities and relative frequencies from the interval $(0, 1)$) is given by
$$w_n(x_1, \dots, x_{t-1}|a_1n, \dots, a_{t-1}n) = C_n \psi(x_1, \dots, x_{t-1}) [x_1^{a_1} \dots x_{t-1}^{a_{t-1}} (1-x_1- \cdots -x_{t-1})^{a_t}]^n, $$ where $\psi$ is the a-priori density with properties analogous to the binomial case, and $C_n$ is determined by the condition
$$\int_{\mathbb{R}^{t-1}} w_n(x_1, \dots, x_{t-1}|a_1 n, \dots, a_t n) dx_1 \cdots dx_{t-1} = 1.$$

In order to apply the multi-dimensional auxiliary theorem, von Mises defined the functions (identical for all $k$)
$$f_k(x_1, \dots, x_{t-1}) = \left\{ \begin{array}{cll} \left(\frac{x_1}{a_1}\right)^{a_1} \cdots \left(\frac{x_{t-1}}{a_{t-1}}\right)^{a_{t-1}} \left(\frac{1-x_1- \cdots - x_{t-1}}{a_t}\right)^{a_t} & \mbox{for} & 0 < \sum_{i=1}^{t-1} x_i < 1\\
0 & \mbox{else.} & \end{array} \right. $$ These functions have a unique maximum in $a:= (a_1, \dots, a_{t-1})$ with $f_k(a) = 1$, and the (automatically) constant matrix $H_n \in \mathbb{R}^{t-1, t-1}$  of the auxiliary theorem is given by
$$h_n^{i,j} = \frac{1}{2a_t} \,\, (i \ne j) \quad h_n^{i,i} = \frac{1}{2} \left( \frac{1}{a_i} + \frac{1}{a_t}\right).$$
Concerning the limit of the distribution function $W_n$ of the posterior von Mises thus arrived at
\begin{multline*} W_n(a_1+\frac{z_1}{\sqrt{n}}, \dots, a_{t-1} + \frac{z_{t-1}}{\sqrt{n}} | a_1n, \dots, a_{t-1} n) =\\
= \sqrt{\frac{\mbox{det} H_n}{(2\pi)^{t-1}}} \int_{-\infty}^{z_1}\!\! \cdots \!\! \int_{-\infty}^{z_{t-1}} \!\! \exp(- x^T H_n x) dx_1 \dots dx_{t-1}.\end{multline*}

It is a notable detail that von Mises's proofs also included corresponding limit theorems for densities. In this respect he was more general than Bernstein.

On the other hand, compared with Bernstein, von Mises held the slightly less general presupposition of constant relative frequencies. For applications, where approximations for large numbers were considered, this difference was not relevant. In this sense \citet{mises:1919a} also gave a rigorous treatment of Bienaymé's result concerning location parameters by means of a corresponding limit theorem. In his 1931 textbook von Mises proved all these theorems directly, without use of auxiliary theorems, but in principle based on the same methods as in 1919. \nocite{mises:1931}

Von Mises's 1919 article on limit theorems is also very notable because it contains a comparison of asymptotic methods as based on direct and on inverse probabilities. The central limit theorem was thus paralleled by the BvM theorem. To my knowledge he was the first author to do so comprehensively and in a systematic way, if still in the tradition of classical problems as they had been pointed out by Laplace.

\section{Neyman's Papers on the Numerical Equality of Frequentist and Bayesian Testing with Large Samples}

In the meantime, Neyman had gone to Poland, where he held different positions at research institutions and at Warsaw University. During a stay in London, supported by a Rockefeller stipend, he met Egon Pearson (1895--1980), and the close collaboration between the two mathematicians began. In 1928 Neyman and Pearson published their first famous paper (in two parts)  in {\it Biometrika}, where, in principle, type I and type II errors were introduced, and also the idea of likelihood ratios as test statistics appeared.

In order to test the hypothesis that a sample $\Sigma$ has been drawn from a ``population'' $\pi$ (from today's point of view $\pi$ would correspond to a particular, possibly multi-dimensional parameter), a test statistic, called ``likelihood'' $$\lambda := \frac{C(\pi)}{C_{\mbox{max}}}$$ was introduced, where, according to the now common use, $C(\pi)$ designated the likelihood function of $\Sigma$ with respect to $\pi$,  and $C_{\mbox{max}}$ the maximum likelihood of $\Sigma$. In a case which Neyman and Pearson [1928] were especially interested in, the ``universe'' of all populations $\Omega$ was corresponding to random variables with discrete values $x_1, \dots, x_k$. In a test on the simple hypothesis that the probability for $x_i$ is equal to $p_i$ ($i=1, \dots, k$, $\sum p_i = 1$) the ``likelihood'' of a sample $\Sigma$ of length $n$  in which the occurrence of $x_i$ is counted $n_i$-times, is
\begin{equation} \label{multiratio} \lambda = \prod_{i=1}^k \left(\frac{p_i}{q_i}\right)^{n_i}, \quad q_i := \frac{n_i}{n}.\end{equation} In a test on the ``composite hypothesis,'' where $\Pi = \omega \subset \Omega$, $\omega$ consisting of more than one element, the likelihood of a sample is defined by
$$\lambda = \frac{\sup_{x \in \omega} C(x)}{C_{\mbox{max}}}.$$

In 1928, Neyman and Pearson also introduced errors of type I and II in a similar way as we know it now. Type I means the probability of rejecting a hypothesis even when it is true. Type II means the probability of accepting a hypothesis even when it is false. They had already the idea that tests based on likelihoods would have a certain optimality concerning type II errors. Therefore the focus was on type I errors, expressed by the probability that under the hypothesis considered the test statistic $\lambda$  becomes $\leq \lambda_0$, where $\lambda_0$ is the likelihood of the drawn sample. In the discrete case under the simple hypothesis $\Pp(x_i)=p_i$  this probability would be
\begin{equation} \label{typeI} P = \Pp_{p_1, \dots, p_k}\left(\left\{(n_1, \dots, n_k) \big| \prod_{i=1}^k \left(\frac{p_i}{q_i}\right)^{n_i} \leq \lambda_0\right\}\right), \quad (\sum n_i = n) \end{equation} according to (\ref{multiratio}). In principle, this probability can be calculated by a multinomial distribution with parameters $(p_1, \dots, p_k)$, but there are unsurmountable difficulties in exactly determining the set of tuples $(n_1, \dots, n_k)$ which meet the condition.

Fortunately, there exists a good approximation of $\lambda$ for large samples. It is relatively easy to show that
\begin{equation} \label{approrat}
\prod_{i=1}^k \left(\frac{p_i}{q_i}\right)^{n_i} \approx \exp\left(- \frac{1}{2} \chi^2\right), \quad \chi^2:= \sum_{i=1}^k \frac{(n_i - p_i n)^2}{p_i n}. \end{equation}
On the other hand, the vector $(n_1, \dots, n_k)$ is approximately normally distributed for large $n$, where the multivariate normal distribution obeys a law proportional to $\exp(-\chi^2/2)$. Thanks to a theorem proved by Karl \citet{pearsonk:1900} one can infer that
\begin{equation} \label{testapprox}
P \approx \frac{\int_{\chi_0}^\infty x^{k-2} e^{-x^2/2} dx}{\int_0^\infty x^{k-2} e^{-x^2/2} dx}, \quad \lambda_0 = \exp(-\chi_0^2/2).
\end{equation}

The significance of the multinomial model consists in the fact that after appropriate grouping also problems concerning continuous variables can be reduced to considering multinomial and, in cases of large samples, multi-variate normal distributions. The $\chi^2$-test on a certain family of distributions is a standard example for a test of this kind, and if the test is on a family of distributions with parameters from a whole range, a composite
hypothesis is considered. In the second part of the 1928 paper this case was comprehensively expounded. Another example would be a test on the mean. When the variable under consideration is split off into discrete values $x_1, \dots, x_k$  with unknown probabilities $p_1, \dots, p_k$, then, the hypothesis $H$---again a composite one---that has to be tested is: \begin{equation} \label{comphyp} H: \mu_0 = x_1p_1+x_2 p_2 + \cdots + x_n p_n.\end{equation} Hypotheses of this kind were treated in \citep{neyman:1930} by means of a-posteriori considerations, as we will see in a moment.

The posterior standpoint was briefly described by Neyman in a paper presented at the International Congress of Mathematicians in Bologna in 1928. This contribution was published only four years later, however, as \citep{neyman:1932}.  By Bernstein's lectures Neyman had learned that an approximate normal distribution also occurs in the context of inverse probabilities, at least in the binomial case. Apparently in 1928 still without knowledge of von Mises 1919 work, he proved a limit theorem for multinomial posteriors, on exactly the same assumptions as von Mises and basically by a similar approach. For corresponding accounts and theorems ``under consideration'' \citet[35]{neyman:1932} referred to articles on Bayes's formula by Karl and Egon Pearson between 1920 and 1924 in {\it Biometrika}, which did not contain rigorous proofs,  unpublished papers by Borel (from his 1926 lectures at the Sorbonne) and by a certain ``Miss A. Miklaszewska in Warsaw'', and finally to Bernstein's lithograph from 1917. Evidently, all these contributions dealt with a less general situation. For a long time this article by Neyman was, as far as I know, the only  ``western'' publication where Bernstein's lithograph was correctly referred to, but it remained widely unnoticed, because, already in 1930, Neyman had published a paper which considerably enlarged the contents of the 1928 presentation.  Concerning von Mises's multinomial limit theorem \citet{neyman:1930} wrote:
\begin{quote} \normalsize
With regard to that theorem I must apologize that I overlooked it for a very long time. Since I proved it independently and was persuaded that it was new, I mentioned it in my paper read before the International Mathematical Congress in Bologna in 1928 without any reference to its original inventor. Unfortunately nobody present noticed that the theorem was not new and its authorship has been wrongly attributed.

Proving the theorem I used a method invented by E. Borel and by S. Bernstein for the proof of a less general theorem.
\end{quote}
At this place, references to Borel's and Bernstein's works were lacking, however.

By means of the asymptotic a-posteriori normality of the multinomial, Neyman \citeyearpar{neyman:1932, neyman:1930} showed an analog to (\ref{testapprox}), namely
\begin{equation} \label{typeIpost}
\Pp(\tilde \chi \geq \tilde \chi_0 | n_1, \dots, n_k) \approx
 \frac{\int_{\tilde \chi_0}^\infty x^{k-2} e^{-x^2/2} dx}{\int_0^\infty x^{k-2} e^{-x^2/2} dx}, \quad {\tilde \chi}^2 := \sum_{i=1}^k \frac{(n_i-p_in)^2}{n_i} .
\end{equation}
Now, in large samples it is to be expected that $n_i \approx np_i$ and therefore $\chi^2 \approx \tilde \chi^2$. This implies that the type I error according to (\ref{typeI}) should be approximately equal to the posterior probability according to (\ref{typeIpost}) for $\lambda_0 = \exp(-\tilde \chi_0^2/2)$. Despite the fundamentally different standpoints in testing by aid of direct or inverse probabilities, numerical results of the respective probabilities are very similar in case of large samples. \citet{neyman:1930} significantly enlarged those considerations to composite hypotheses of the type (\ref{comphyp}). This was Neyman's own ``particular attack on the Bayes front.''

\section{Le Cam's General Theorem}

After having graduated from Sorbonne in 1947, Lucien Le Cam (1924--2000) until 1950 worked as a statistician, but then went to the University of California at Berkeley and there continued his studies.
His doctoral thesis, supervised by Neyman and published in 1953, especially focused on a generalization of the concept of efficiency of parameter estimators by considering ``Bayes estimators.'' (Le Cam still used there the denotation ``estimate'' both for estimators as r.v's and for their realizations.) In order to discuss optimality he introduced---at first not closer defined---gain functions  $W_k(T_k(Z_k), \theta)$, where $T_k$ is an estimator for the parameter $\theta \in \Theta \subset \mathbb{R}^r$ defined on the sample vector $Z_k$ of dimension $k$ consisting of  i.i.d. coordinates $X_1, \dots, X_k$ with a common probability density $f_\theta$ with respect to a measure $\nu$. (Here and in the following, technical details on measurability and similar issues are omitted; in order to unify the exposition only real-valued sample elements are considered.) Utility functions $R(T_k, \theta)$ are then given by
$$R(T_k, \theta) := \Ee(W_k(T_k(Z_k), \theta)),$$ and criteria for goodness of estimators refer to the expression
$$J(T_k):=\int_\Theta R(T_k, \theta) \lambda(\theta) d \theta,$$ where $\lambda(\theta)$ is a non-negative weight function. The integral $J(T_k)$ can then be considered as an average utility. If  for $z_k=(x_1, \dots, x_k)$ we set
$$p_k(z_k, \theta) = \prod_{i=1}^k f_\theta (x_i),$$ such that $p_k(\cdot, \theta)$ is a density with respect to the product measure $\nu_k := \nu^k$, we obtain
$$\int_\Theta R(T_k, \theta) \lambda(\theta) d \theta = \int_\Theta \int_{\mathbb{R}^k} W_k(T_k(z_k), \theta) p_k(z_k, \theta) \lambda(\theta) d \nu_k(z_k) d \theta.$$ If we interpret $\lambda(\theta)$ as an a-priori density for $\theta$, then with the abbreviation $p_k(z_k) := \int_\Theta p_k(z_k, \theta) \lambda(\theta) d \theta$ it follows
\begin{equation} \label{bayesest} J(T_k)=\int_\Theta R(T_k, \theta) \lambda(\theta) d \theta =  \int_\Theta \int_{\mathbb{R}^k} W_k(T_k(z_k), \theta) p_k(z_k) p_k(\theta | z_k) d \nu_k(z_k) d\theta.\end{equation} This is the reason why estimators with certain optimality properties depending on $J(T_k)$ are called Bayes estimators. More precisely: If for the estimator $\mu_k$ the relation $J(\mu_k) = \sup_{T_k \in S_k} J(T_k)$ is valid, where $S_k$ is the set of all estimators for which $J(T_k)$ in the sense of an integral exists, $\mu_k$ is named $k$-Bayes estimate; and if for a sequence of estimators $(\mu_k)$ the inequality $J(\mu_k) \geq \sup_{T_k \in S_k}J(T_k) - \varepsilon_k$ holds, such that $\epsilon_k > 0$ tends to 0 for $k \to \infty$, the sequence $(\mu_k)$ is called an asymptotic Bayes estimator.

If the inverse density $p_k(\theta | z_k)$ tends to a normal density in an appropriate way, and after choosing an appropriate utility function, (\ref{bayesest}) may lead to an asymptotic optimality of $J(T_k)$, $T_k$ being a maximum likelihood estimator. In order to tackle this problem, Le Cam proved a general BvM theorem in the sense of our sect. 1.

Under certain conditions there exists an open subset $C \subset \Theta$ and a sequence $(\mu_k)$ of uniquely determined maximum likelihood estimators with $\Pp(\mu_k(z_k) \to \theta|\theta) = 1$ for all $\theta \in C$, as Le Cam showed first. Then, if the a priori density $\lambda$ is continuous and positive in $\theta_0 \in C$, and if the density $f(x, \theta)$ meets certain regularity conditions, in particular continuous second order derivatives with respect to $\theta$, such that for the $r \times r$ matrix $B(x, \theta) := - \frac{\partial}{\partial \theta_i \partial \theta_j} \log f(x, \theta)$ the matrix $\Gamma(\theta) = \int_{\mathbb{R}} B(x, \theta) f(x, \theta) d x$ is positive definite,  then the limit relation
\begin{equation} \label{lecamtheorem} \Pp\left(\lim_{k \to \infty} \int_{\mathbb{R}^r} \left|\frac{1}{\sqrt{k^r}} p_k(\mu_k-\frac{t}{\sqrt{k}}|z_k)-\frac{(\mbox{det} \Gamma(\theta_0))^{1/2}}{(2\pi)^{r/2}}\exp(- t^T \Gamma(\theta_0) t)\right| d t = 0 \big| \,\, \theta_0\right) = 1 \end{equation} is valid. In this limit relation it is supposed that $p_k(\mu_k-\frac{t}{\sqrt{k}}|z_k) = 0$ if $\mu_k-\frac{t}{\sqrt{k}} \not \in C$. Important for applications was the (very strong) property that the mode of convergence of distributions was in (total) variation, as indicated here through the integral $\int_{\mathbb{R}^r} | \cdot | dt$.

Especially useful to the application of this limit theorem to asymptotic optimality, particularly of maximum likelihood estimators, were gain functions of the type
$$W_k(T_k, \theta) \propto \exp((T_k-\theta)^T D (T_k-\theta)), $$ where $D$ is a positive definite $r \times r$-matrix.

In the 1953 paper and in later works (e.g., 1986, 1990/2000) Le Cam acknowledged the merits of Laplace, Bernstein, and von Mises regarding the asymptotic normality of the posterior. In the last of these works, written jointly with Grace Lo Yang, notions like ``Bernstein-von Mises phenomenon'' or ``Bernstein-von Mises result'' were explicitly introduced. At the same time Le Cam was eager to point out the significance of his own, by far more general, theorem. This effort was, however, combined with some historical inaccuracies or even inconsistencies. It is not the goal of the present paper to pedantically list the problematic points in detail. The interested reader may compare Le Cam's historical notes in the above-mentioned works with our sects. 2--4. Evidently, these notes have to be understood from the point of view of a working mathematician: Le Cam wanted to emphasize the general significance of his own theorem compared with previous attempts, and, of course, there is no doubt about this significance. Still, his assessment of Bernstein's account as a work which has ``little to do with that of Laplace'' \citeyearpar[261]{lecam:2000} seems a little strange. With certainty Le Cam owed the information about Bernstein's proof to his advisor Neyman, but it may be that this information was only superficial. This impression is evidenced by the fact that references to Bernstein in Le Cam's work read always ``Theory of probability, 1917 (Russian).'' These incomplete data have caused some confusion, at least among readers outside Russia, because a 1917 book with this title could not be found in library catalogues. Therefore it remained unclear for a long time what character this ominous 1917 publication by Bernstein really had.

\section{Conclusion}

As we have seen, the basic idea behind proving the asymptotic normality of posterior distributions is due to Laplace, and from the very beginning this result served for establishing favorable properties of estimators and for testing hypotheses. Within these problem fields, Laplace and his direct successors changed in a rather pragmatic way between a-priori and a-posteriori considerations. Beside this, for Laplace the increasing concentration of the posterior in a small neighborhood around the relative frequency served as a model for  the epistemic process that leads to the determination of probabilities by counting frequencies. Bernstein and von Mises possess the merit of rigorous proofs (according to now usual standards) in particular however very important situations. Von Mises was apparently the first to give a systematic exposition of asymptotic methods for direct and inverse probabilities, still in the context of traditional probability theory, whereas Neyman applied the multinomial version of the BvM for a discussion of particular questions regarding his (and Egon Pearson's) new approach to hypothesis testing. Le Cam considerably enlarged and generalized the assertion of the limit theorem and used it to the assessment of goodness of maximum likelihood estimators.

Altogether, the name ``Bernstein-von Mises theorem'' seems to be not inappropriate, but a designation like ``Laplace-Le Cam theorem'' would be at least equally fitting.

\bigskip

\noindent {\it Acknowledgment}: I thank an unknown colleague from the University of Kharkiv for providing me with a copy of Bernstein's 1917 lithograph. To Reinhard Siegmund-Schultze (Kristiansand, Norway) I owe valuable information about von Mises's life and work, and he drew my attention to the quotation about Bernstein in Reid's book. Richard Pulskamp (Cincinnati, OH) gave very useful suggestions which helped to improve the text.

\end{document}